\documentclass[11pt]{amsart}
\usepackage{geometry}                
\usepackage{graphicx}
\usepackage{amssymb}
\usepackage{epstopdf}
\newtheorem{thm}{Theorem}

\newtheorem{prop}[thm]{Proposition}

\newtheorem{lemma}[thm]{Lemma}

\newtheorem{corollary}[thm]{Corollary}

\theoremstyle{definition}

\newtheorem{defn}[thm]{Definition}

\theoremstyle{remark}

\newtheorem{example}{Example}
\newtheorem{rem}{Remark}

\let\oldmarginpar\marginpar
\renewcommand\marginpar[1]{\-\oldmarginpar[\raggedleft\footnotesize #1]%
{\raggedright\footnotesize #1}}

\newcommand{\g}{{\mathfrak{g}}}                                            

\newcommand{\C}{\mathbb C}
\newcommand{\End}{{\mathrm{End}}}

\newcommand{\ot}{\otimes}

\newcommand{\ZZ}{\mathbb Z}

\newcommand{\cR}{\mathcal{R}}
\newcommand{\ad}{\textrm{ad}}

\renewcommand{\gg}{\mathfrak{g}}

\newcommand{\eps}{\epsilon}
\newcommand{\NN}{\mathbb{N}}
\renewcommand{\H}{\mathcal{H}}
\newcommand{\HH}{\mathcal{H}\!\!\!\mathcal{H}}
\newcommand{\cC}{\mathcal{C}}

\newcommand{\CC}{\mathbb{C}}
\newcommand{\id}{\textrm{id}}

\title{Quantum D-modules, elliptic braid groups and double affine Hecke algebras}
\author{David Jordan}
\begin{document}
\maketitle

\begin{abstract}
We build representations of the elliptic braid group from the data of a quantum $D$-module $M$ over a ribbon Hopf algebra $U$.  The construction is modelled on, and generalizes, similar constructions by Lyubashenko and Majid \cite{Ly}, \cite{LyMa}, and also certain geometric constructions of Calaque, Enriquez, and Etingof \cite{CEE} concerning trigonometric Cherednik algebras.  In this context, the former construction is the special case where $M$ is the basic representation, while the latter construction can be recovered as a quasi-classical limit of $U=U_t(sl_N)$, as $t\to 1$.  In the latter case, we produce representations of the double affine Hecke algebra of type $A_{n-1}$, for each $n$.
\end{abstract}

\section{Introduction}
Let $G=SL_N(\C)$, and $V$ be the vector representation of $G$. In
\cite{AS}, Arakawa and Suzuki construct a functor $F_n$ from the
category of $U(\g)$-bimodules to the category of representations
of the degenerate affine Hecke algebra (AHA) $\H_n^{deg}$. Namely,
as a vector space, $F_n(M)=(V^{\otimes n}\otimes M)^\g$ (where
$\g={\rm Lie}(G)$ acts on $M$ by the adjoint action), and the
generators of $\H_n^{deg}$ act on $F_n(M)$ by certain explicit
formulas.

In the paper \cite{CEE}, Calaque, Enriquez, and Etingof
extended this construction to the double
affine case. Namely, they upgraded the Arakawa-Suzuki functor to a
functor $F_n$ from the category of $D_G$-modules to the category
of representations of the degenerate (i.e., trigonometric) double
affine Hecke algebra (DAHA) $\HH_n^{deg}(k)$, $k=N/n$. They 
also considered the rational degeneration of this construction, in
which one uses $D_\g$-modules instead of $D_G$-modules, and the
rational DAHA $\HH_n^{rat}(k)$ instead of the trigonometric one.

The goal of the present paper is to generalize both of these
constructions to the case of quantum groups and nondegenerate AHA
and DAHA.  Namely, let $U$ be a ribbon Hopf algebra with enough
finite dimensional representations (i.e. the common annihilator
of all finite dimensional representations is zero), and let $V$ be
an irreducible finite dimensional representation of $U$. Then we
define a functor $F_{n,V}$ from the category of $U$-bimodules to
the category of representations of the affine braid group
$B_n^{tr}$, given by the formula $F_{n,V}(M)=(V^{\otimes
n}\otimes M)^{\rm inv}$, where the invariants are taken with
respect to the adjoint action of $U$ on $M$, and the action of
the generators of the affine braid group is defined using R-matrices. These
representations are similar to those considered in
\cite{EtGe}. If $V$ satisfies the Hecke condition
(i.e. the braiding on $V\otimes V$ satisfies the equation
$(\beta-q^{-1}t)(\beta+q^{-1}t^{-1})=0$), then $F_{n,V}(M)$
descends to a representation of the AHA $\H_n(t)$.  More
interestingly, we upgrade the functor $F_{n,V}$ to a functor from
the category of $D_U$-modules ($D$-modules on the quantum group
corresponding to $U$) to the category of representations of the
elliptic braid group $B_n^{ell}$, which in the Hecke case lands
in the category of representations of the DAHA $\HH_n(q,t)$. If
$U=U_t(sl_N)$, $V$ the vector representation, and $t=q^{nk}$,
then in the quasiclassical limit $q\to 1$ we recover the functors
from \cite{AS} and \cite{CEE}, respectively.

Our construction is also a generalization of the work of
Lyubashenko and Majid \cite{Ly}, \cite{LyMa}, where an action of the elliptic
braid group is obtained on $(V^{\otimes n}\otimes A)^{\rm inv}$,
where $A$ is the dual Hopf algebra of $U$. Indeed, $A$ is the
most basic example of a $D_U$-module (the module of functions on
the quantum group).

Finally, we would like to discuss the connection of our paper
with the work of Varagnolo and Vasserot, \cite{VV}. This
connection is the quantum counterpart of the connection between
the results of \cite{CEE} and those of \cite{GG}. Namely,
consider the setting of the present paper with $U=U_t(sl_N)$,
$V=\Bbb C^N$, and let $e_n$ be the Young symmetrizer of the
finite Hecke algebra contained in $\HH_n(q,t)$. Then the spherical
subalgebra $e_n\HH_n(q,t)e_n$ acts in the space
$e_nF_{n,V}(M)=(S^n_tV\otimes M)^{\rm inv}$ where $S^n_tV$ is the
quantum symmetric power. Also, recall that Varagnolo and Vasserot
define a functor $\Phi_{N,k}$ from the category of D-modules on
the product $G_t\times \Bbb P_t^{N-1} $ of the quantum group
$G_t$ with the quantum projective space $\Bbb P_t^{N-1}$ twisted
by the $k$-th power of the line bundle $O(1)$, to the category of
representations of the spherical DAHA $e_N \HH_N(q,q^k)e_N$.  It
is easy to see that the space of global sections of $O(1)^n=O(n)$
over the quantum projective space is $S^n_tV$, so $\Phi_{N,n}(M)=
e_n F_{n,V}(M)$, and thus through the combination of the
constructions of the present paper and \cite{VV}, two different
algebras $e_n \HH_n(q,q^N)e_n$ and $e_N \HH_N(q,q^n)e_N$ get to act
on the same vector space. In a future paper, we plan to show that
the images of these actions in the endomorphism algebra of this
space are the same; in the rational degeneration, this was proved
in \cite{CEE}. In particular, this will imply that if $N$ divides
$n$, then $e_N \HH_N(q,q^n)e_N$ is a quotient of $e_n
\HH_n(q,q^N)e_n$. Moreover, recall from \cite{VV2} that the
algebra $e_n \HH_n(q,t)e_n$ admits an ``analytic continuation'' with
respect to $n$, which yields the quantum toroidal algebra
$Q_\lambda(q,t)$ of type $gl(1)$, which projects to $e_n \HH_n(q,t)e_n$ when
$\lambda=n$ is a positive integer. The above statements should
follow from the existence of an isomorphism
$Q_\lambda(q,q^\mu)\to Q_\mu(q,q^\lambda)$.

In subsequent papers, we plan to study the
representation-theoretic properties of the functor $F_{n,V}$,
i.e. what it does to particular $U$-bimodules and
$D_U$-modules; in particular, it would be interesting 
to consider the case of roots of unity. 
Also, our construction paves way for a number of
further generalizations, which we plan to explore in the
future. One of them is the quantum generalization of the paper
\cite{EFM}, which generalizes the construction of \cite{CEE} from
the type A case to the type BC case: it defines a functor from
twisted D-modules on the symmetric space $GL_N/GL_p\times GL_q$
($p+q=N$) to representations of degenerate DAHA of type
$BC_n$. This generalization will involve quantum symmetric spaces
and the Sahi-Stokman $BC_n$ DAHA, and will be related to the
paper \cite{OS} in the way similar to the relation
between the present paper and the construction of \cite{VV},
explained above. Another interesting direction is the
generalization to the case of an arbitrary ribbon category
$\mathcal C$, which is not necessarily the category of
representations of a ribbon Hopf algebra. For example, if this
category is semisimple, the role of the dual $A$ of $U$ will be
played by $\oplus_X X\otimes X^*$, where $X$ runs over all simple
objects of $\mathcal C$. In this case, if $M=A$, our construction
would recover the natural elliptic braid group action on the
genus 1 modular functor, described in the book \cite{BK}. We
note that this elliptic braid group action comes with a
compatible action of the modular group $SL_2(\Bbb Z)$, and we
expect that under some conditions on the D-module $M$, such an
action will exist on $F_{n,V}(M)$; in the Hecke case this will
recover the difference Fourier transform of Cherednik
\cite{Che}. Finally, we would like to use the approach of this paper to
construct and study quantization of multiplicative quiver varieties of
Crawley-Boevey and Shaw \cite{CBS}.

The contents of this article are laid out as follows.  Section 2 consists of preliminaries.  In Subsection 2.1, we recall the definitions of the elliptic braid group and a certain quotient, called the double affine Hecke Algebra.  In Subsection 2.2, we recall the notion of twisting of the comultiplication of a quasi-triangular Hopf algebra $U$.  This subsection is somewhat technical, and may be skipped on a first read.  In Subsection 2.3, we recall the construction of the Reflection Equation algebra, which in a certain sense generalizes the algebra of functions on an algebraic group to the braided setting.  To the extent possible, we give key definitions and propositions in braided-categorical terms.  In Subsection 2.4, we recall the Heisenberg double construction, and its relation to differential operators on an algebraic group.  In Subsection 2.5, we recall how to apply this in the non-commutative (quantum, braided) context.  In particular, here is where we give the definition of quantum $D$-modules we will use in this article.  In Subsection 2.6 we recall the left, right, and adjoint actions of a Lie algebra $\gg$ on $D_G$-modules, and generalize this to the quantum setting.

In Section 3, we state without proof our main results; in particular, we assert the existence of a certain family of functors from the category of quantum $D$-modules to the category of representations of the elliptic braid groups.  Section 4 comprises our primary new contribution to the subject, wherein we construct the functors $F_{n,V}$ asserted in Section 3.  In Section 5, we explain that the element $\tilde{Y}$ of Section 4 and a related element $\tilde{X}$ act as scalars when $V$ is irreducible.  In Section 5, we consider the case $U=U_t$, and we show that we can recover the geometric constructions of \cite{CEE} as a trigonometric degeneration of our constructions.  

\subsection{Acknowledgments}
I would like to heartily thank Pavel Etingof for explaining the construction \cite{CEE} in the degenerate case, for his considerable help with the present construction, and for his contribution to the introduction.  I would also like to thank Kobi Kremnizer for countless helpful conversations where I learned about classical and quantum $D$-modules.  This project would have proceeded nowhere without their guidance.  The author's work was partially supported by the NSF grant DMS-0504847.

\section{Preliminaries}
\subsection{The Elliptic Braid Group and the DAHA}
In this section we define the elliptic braid group following Birman and Scott \cite{Bir} and \cite{Sc}, and a particular quotient of its group algebra, called the double affine Hecke algebra.
\begin{defn}  \label{BnEllRelns}The elliptic braid group, $B_n^{Ell}$, is the fundamental group of the configuration space of $n$ points on the torus.
It is generated by
\begin{itemize}\addtolength{\itemsep}{-0.6\baselineskip}
\item the commuting elements $X_1,\ldots,X_n$\\
\item the commuting elements  $Y_1,\ldots,Y_n$\\
\item and the braid group of the plane,
\begin{eqnarray*}B_n=\left<T_1,\ldots, T_{n-1} \phantom{--}\Big| \phantom{--} \begin{array}{l}T_iT_{i+1}T_i=T_{i+1}T_iT_{i+1} \forall i,\\T_iT_j=T_jT_i, |i-j|\geq 2\end{array}\right>.
\end{eqnarray*}
\end{itemize}
The cross relations are:
\begin{itemize}\addtolength{\itemsep}{-0.6\baselineskip}
\item $T_iX_iT_i=X_{i+1},$\\
\item $T_iY_{i}T_i=Y_{i+1},$\\
\item $X_1Y_2=Y_2X_1T_1^2$\\ 
\item $\tilde{Y} X_i = X_i \tilde{Y}$, where $\tilde{Y}=\prod_jY_j$.
\end{itemize}
\end{defn}

\begin{figure}[htb]
\centerline{\includegraphics[height=2.5cm]{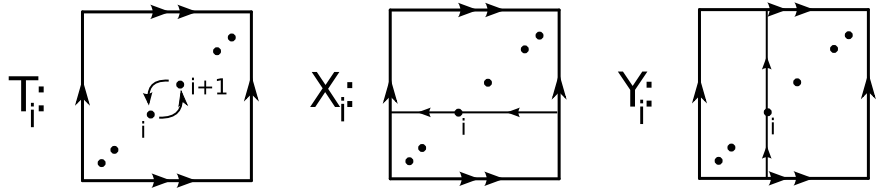}}
\caption{Generators for the elliptic braid group $B_n^{Ell}$.}\label{BnEll}\end{figure}

Under the usual realization of the torus $T^2$ as the unit square with opposite sides glued, we choose as a basepoint the configuration with all the marked points along the diagonal.   The subgroup $B_n$ is then identified with those braids which stay away from the sides, while the $X_i$ and $Y_i$ correspond to the horizontal and vertical global cycles, respectively.  This is depicted in Figure \ref{BnEll}.  \footnote{The reader should note that these are slightly non-standard generators (e.g. \cite{Che}).  In particular, the $Y_i$ here are the inverses of the $Y_i$ there.}

\begin{defn} \label{DAHAdefn} Cherednik's double affine Hecke algebra (DAHA) $\HH(q,t)$ is the quotient of the group algebra $\CC B_n^{Ell}$ of the elliptic braid group by the additional relations
$$(T_i-q^{-1}t)(T_i+q^{-1}t^{-1})=0,$$
where $q,t$ are complex parameters.  Note that when $q,t=1$, we have an isomorphism $\HH(1,1)\cong\CC[S_n\ltimes \ZZ^{2n}]$.
\end{defn}

\begin{rem}  Let us discuss the precise relation between 
our DAHA $\HH_n(q,t)$ and the one appearing in \cite{Che}. 
One can define a 3-parameter DAHA, $\HH_n(q,t,p)$ by replacing the last 
relation of $\HH_n(q,t)$ with $\tilde YX_i=pX_i\tilde Y$. Then our algebra 
is  $\HH_n(q,t,1)$ and Cherednik's is $\HH_n(1,t,p)$.  
\end{rem}

\subsection{Twistings on a Quasi-triangular Hopf algebra}

Throughout this article, $U$ denotes a quasi-triangular Hopf algebra with universal $R$-matrix $\cR=\sum_kr^+_k\ot r^-_k$. We will denote by $\Delta, \eps, S$ the comultiplication co-unit, and antipode in $U$.  Let $(\cC,\ot,\sigma)$ denote the braided tensor category of its left modules.  Here $\sigma$ is given by $\tau\circ \cR$, where $\tau$ is the flip of tensor factors.  

We will need the Hopf algebras $U^o$, $U^e$, and $U^{[2]}$ as they appear in \cite{VV}.  For clarity's sake, we adopt the same notations, and recall them here.

\begin{defn} Let $U^o$ denote the algebra $U$, with opposite co-multiplication, $\Delta^o(x)=\tau\circ\Delta(x)$.
\end{defn}

\begin{defn} Let $U^e$ denote the algebra $U^o\ot U$, with coordinate-wise multiplication, co-multiplication, and antipode.
\end{defn} 

\begin{defn} Let $U^{[2]}$ denote the algebra $U\ot U$ with coordinate-wise multiplication, but with comultiplication given by $\tilde{\Delta}(x\ot y)=\cR_{23}^{-1}\tau_{23}(\Delta(x)\ot\Delta(y))\cR_{23}$, and with antipode $\tilde{S}(x\ot y)=\cR_{21}(S(x)\ot S(y))\cR_{21}^{-1}$.\end{defn}

In fact, $U^e$ and $U^{[2]}$ are related by a 2-cocycle, which induces an equivalence on their tensor categories of modules.  Let us recall these constructions:

\begin{defn} Let $H$ be a Hopf algebra.  A normal left 2-cocycle on $H$ is an invertible element $c\in H\ot H$ such that
$$(\eps\ot \id)(c)=(\id\ot \eps)(c)=1\ot 1, \textrm{ and } (\Delta\ot \id)(c)(c\ot 1)=(\id\ot\Delta)(c)(1\ot c).$$
\end{defn}

A 2-cocycle is sometimes called a twist.  Given a 2-cocycle $c$ on $H$, we can define the twisted Hopf algebra $H_c$ to be the algebra $H$ with twisted comultiplication $\Delta_c(h)=c^{-1}\Delta(h)c$ and antipode $S_c(h)=QS(h)Q^{-1}$, where $Q=\mu\circ(\id\ot S)(c)$. We have the following standard proposition:

\begin{prop} $c$ induces a tensor equivalence $H$-mod$\to H_c$-mod.\end{prop}

It is now straightforward to check that $c=\cR_{13}\cR_{23}$ is a 2-cocycle for $U^e$, and that $U^{[2]}=U^e_c$.

\begin{defn} Let $H$ be a Hopf algebra.  An $H$-equivariant algebra $A$ is an algebra $A$ with an $H$-action on $A$ s.t. the product $\mu:A\ot A\to A$ is a map of $H$-modules.
\end{defn}

\begin{rem} In other words, $A$ is an algebra in the category of $H$-modules.\end{rem}

Given an $H$-equivariant algebra $A$, we can define an $H_c$-equivariant algebra $A_c$ as the same underlying $H$-module, with multiplication given by \mbox{$\mu_c(a\ot b) = \mu(c(a\ot b))$}.  We call $A_c$ the $H_c$-algebra \emph{equivalent} to the $H$-algebra $A$.

For a $(U,U)$-bimodule, we denote the left and right actions by $\rhd, \lhd$ (e.g $x\rhd m \lhd y$).  A $(U,U)$ bimodule is the same thing as a $U^e$-module under the identification $(a\ot b)(v) = b\rhd v \lhd S(a)$.  The co-multiplication on $U$ gives the algebra maps $\Delta: U\to U^{e}$ (resp. $U^{[2]}$).  Thus for a $U^e$-  (resp. $U^{[2]}$-) module $V$, we have an action of $U$, denoted ``$\ad$" given by $(\ad x)(v)=\Delta(x)(v)$ (we use the symbol ``$\ad$" in both contexts).
\subsection{The Reflection Equation Algebra}

In this section, we will recall the so-called reflection equation algebra $A$ associated to the quasi-triangular Hopf algebra $U$.  In the case that $U=U(\gg)$ is the universal enveloping algebra of a Lie algebra $\gg$ of an algebraic group $G$, $A$ will be the algebra functions on $G$.  When $U=U_t(\gg)$ is the quantum group associated to the Lie algebra $\gg$\footnote{Note that $U_t$ is not quasi-triangular, since the $R$-matrix lies in a completion of $U_t\ot U_t$.  However, in all our constructions, we always apply one of the components of $\cR$ to a finite dimensional module, so its action is well-defined.}, $A$ will be a braided version of the algebra of functions on $G$, distinct from the dual quantum group $O_q$.  Majid called $A$ the ``the braided Hopf algebra associated to $O_q$".  The primary advantage of $A$ from our perspective is that there is an adjoint action of $U$ on $A$, for which the algebra structure on $A$ is equivariant, and which does not exist for the usual dual quantum group $O_q$.  This equivariance property was first observed and explained by Majid \cite{Maj}, who proposed the reflection equation algebra as a preferable replacement for $O_q$ in the context of braided differential geometry, and showed that it was a braided-commutative braided-Hopf algebra in the category of $U$-modules.

The results here are all standard, and can be found in one form or another in many sources, e.g. \cite{Maj} or \cite{KlSch}.  We include them here for completeness, and to establish the diagrammatical notation which will appear in later sections.  Also, the reflection equation, Proposition \ref{refleqn}, is usually stated for the defining modules for the FRT bialgebra, but we will need it for arbitrary modules, and so we give a diagrammatical proof.

\begin{defn} Let $F$ denote the restricted dual Hopf algebra of $U$ relative to the tensor category of its finite dimensional representations.  It is the sum of all finite dimensional $U^e$-submodules of $U^*$, and is spanned by functionals $c_{f,v}$, defined by $c_{f,v}(x)=f(xv)$, where $f\in V^*, v\in V$, and $V$ is a finite dimensional $U$-module.
\end{defn}
One can easily check that $c_{f,v}c_{g,w}=c_{f\ot g,v\ot w}$, so that $F$ is a sub-algebra.  The $U^e$-action $(x\ot y)c_{f,v}=c_{xf,yv}$ makes $F$ into a $U^e$-equivariant algebra.  This corresponds to the natural $(U,U)$-bimodule structure on $U^*$ given by $(x\rhd \phi \lhd y)(h)=\phi(S(x)hy)$.

In fact, $F$ is a Hopf algebra with co-product $\Delta(c_{f,v})=c_{f,e_i}\ot c_{e_i^*,v}$, where $e_i$ is a basis for $V$ and $e_i^*$ is a dual basis.  The antipode on $F$ is the adjoint to the antipode on $U$, defined by $(S c_{f,v})(u)=c_{f,v}(Su)$.

\begin{prop}\label{auto} Let $\phi:V\to V$ be a $U$-module map, and let $\phi^*$ denote the adjoint map.  Then $c_{f,\phi v} = c_{\phi^*f,v}$.
\end{prop}
\begin{proof}
For $x\in U$, we compute,
\begin{displaymath}
c_{f,\phi v}(x) = f(x\phi(v))= f(\phi(xv))=(\phi^*f)(xv)=c_{\phi^*f,v}(x)
\end{displaymath}
\end{proof}

\begin{defn} A dual pairing of two Hopf algebras $H$ and $K$ is a map $\kappa:H\ot K\to \CC$ s.t. for all $h,h'\in H$, $k,k'\in K$, we have
\begin{eqnarray*}
\kappa(\Delta_H(h),k\ot k') &=& \kappa(h,kk')\\
\kappa(hh',k)&=&\kappa(h\ot h',\Delta_K(k))\\
\kappa(h,1_K)&=&\eps_H(h)\\
\kappa(1_H,k) &=&\eps_K(k)\\
\kappa(S(h),k)&=&\kappa(h,S(k))
\end{eqnarray*}
The pairing is called non-degenerate if its left and right kernels are zero.
\end{defn}
\begin{defn}\label{nondeg} We say that $U$ \emph{has enough finite-dimensional representations} if the common annihilator of all finite-dimensional representations is zero.
\end{defn}

\begin{rem} The natural pairing of vector spaces $\kappa:U\ot U^*\to\CC$ restricts to a dual pairing of Hopf algebras $U$ and $F$.  Since $F$ was defined as a subalgebra of $U^*$, the right kernel of $\kappa$ is automatically zero.  The left kernel of $\kappa$ is the set of $x\in U$ s.t. $c_{f,v}(u)=0$ for all $f\in V^*, v\in V$, $V$ a f.d. $U$-module.  Thus the natural pairing is non-degenerate if, and only if, $U$ has enough finite-dimensional representations.  This is a mild condition which is certainly satisfied for any $U_t$, and which can be forced in general by quotienting $U$ by the left kernel of $\kappa$, which will be a Hopf ideal.  In this case, the natural pairing between $F$ and $U$ is nondegenerate, and we get $$\textrm{Finite dimensional $U$-modules} \cong \textrm{Finite dimensional $F$-comodules}$$  We will assume from now on that $U$ has enough finite dimensional representations.
\end{rem}

\begin{defn}\label{Adefn} We denote by $A$ the $U^{[2]}$ algebra equivalent to the $U^e$-algebra $F$ via the cocycle $c=R_{13}R_{23}$.  It has the same co-multiplication as $F$, but its multiplication $\mu'$ is related to that of $F$ by the formula:
\begin{eqnarray}
\mu'(f\ot g) = \sum_k \mu(\ad r^+_k(f)\ot (g\lhd S(r^-_k)))\label{multdefn}
\end{eqnarray}
\end{defn}

\begin{rem}\label{inv-desc} For any $U$-module $V$, we have morphisms of $U$-modules
\begin{eqnarray*}
c_{-,-}:V^*\ot V\to A,
f\ot v \mapsto c_{f,v}.
\end{eqnarray*}
Many quantities are most easily computed in the pre-image $V^*\ot V$, applying Proposition \ref{auto} as needed.  In particular, the structure maps $(\mu, \Delta, S)$ for the braided-Hopf algebra $A$, as well as the proofs of Propositions \ref{refleqn} and \ref{mainformula} will be given in this way.
\end{rem}

\begin{figure}[htb]
\centerline{\includegraphics[height=2.5cm]{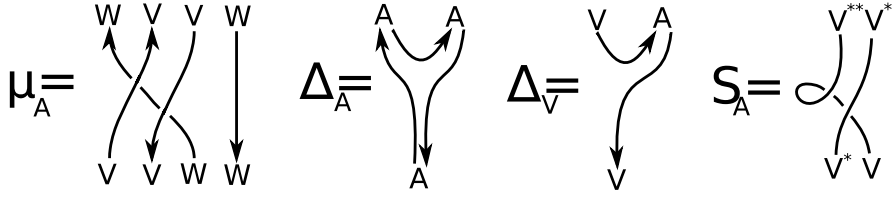}}
\caption{The multiplication and co-multiplication in $A$, the co-action on a $U$-module $V$, and the antipode in $A$.  Diagrams are to be read from bottom up.  The identity morphism $id_V:V\to V$ is denoted with a downward flowing arrow, while $id_{V^*}:V^*\to V^*$ is denoted with an upward flowing arrow.  To avoid confusion, we omit arrows and explicitly label duals in the description of the antipode.  An excellent reference for the diagrammatical calculus in braided tensor categories is \cite{Ka}.}\label{coactmult}\end{figure}

\begin{rem}Let $V$ be a $U$-module.  By duality, $V$ is an $A$-comodule, so we have a canonical element $L \in End_\CC(V)\ot A$ given by $\Delta(v) = L_1(v)\ot L_2$.  This coaction may be described invariantly under the above identification as simply $\Delta: v\mapsto  \sum_ie_i\ot e_i^*v$.  Here $\sum_ie_i\ot e_i^*$ has an invariant description as the image of $1\in\CC$ under the coevaluation map $coev: \CC\to V\ot V^*$. Similarly, formula (\ref{multdefn}) has an invariant description in terms of the braiding.  These are depicted in figure \ref{coactmult}.\end{rem}

\begin{prop}\label{refleqn} \cite{Maj}   The element $L$ defined above satisfies the reflection equation $$L_{01}R_{12}L_{02}R_{21}= R_{12}L_{02}R_{21}L_{01}$$ in the space of endomorphisms of the module \mbox{$V\ot V\ot A$}.  (the tensor indices run from right to left, with $A$ in index 0)\end{prop}
\begin{proof}
The proof is presented in figure \ref{refleqndiag} in the braided tensor category $U$-mod.  The LHS diagram is $L_{01}R_{12}L_{02}R_{21}$, and the RHS is $R_{12}L_{02}R_{21}L_{01}$.  Note that on the RHS we have applied Proposition \ref{auto} to the automorphism $\sigma_{VV}$ of $V\ot V$.
\begin{figure}[htb]
\centerline{\includegraphics[height=4cm]{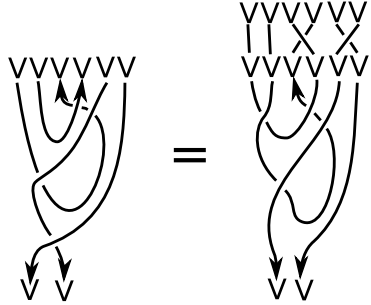}}
\caption{A braided-categorical proof of the ``Reflection Equations".}\label{refleqndiag}
\end{figure}
\end{proof}

\subsection{The Heisenberg double and differential operators}
\begin{defn}  Let $F$ and $U$ be dually paired Hopf algebras.  Then $F$ is an $U^o$-module under the action $x\ot y \mapsto \kappa(S(x),y_1) y_2$, which makes $F$ into a $U^o$-equivariant algebra.  The Heisenberg double $H(F,U)$ is the semi-direct product of $F$ and $U^o$ under this action.  As a vector space it is $F\ot U$, with subalgebras $F\cong F\ot 1$, $U \cong 1\ot U$, and cross relations
$$x f=(x_1\rhd f) x_2=\kappa(S(x_1),f_1)f_2x_2\textrm{, for all }x\in U, f\in F$$
\end{defn}

\begin{example} $D$-modules on an affine algebraic group.
The algebra of differential operators on an affine algebraic group has a particularly satisfactory description in terms of the Heisenberg double.  Let $U=U(\gg)$ denote the universal enveloping algebra of the Lie algebra $\gg$ of an algebraic group $G$, and let $O=O(G)$ denote the algebra of polynomial functions on $G$.

The Lie algebra $\gg$ is constructed as the sub-Lie algebra of vector fields on $G$ which are left-translation invariant.  Thus we have a pairing $\kappa: U\ot O\to \CC$ given by $X\ot f\mapsto X(f)_{\id}$.  This is a dual pairing of Hopf algebras, and so we may construct the Heisenberg double $H(O,U)$.  For instance, if the group $G$ is the affine plane $\CC^n$, one finds $O=\CC[x_1,\ldots,x_n]$, $U=\CC[\partial_1,\ldots,\partial_n]$, and $H(O,U)$ is the $n$th Weyl algebra $W_n$.  More generally, we have the following standard proposition (see, e.g. \cite{KlSch}):
\begin{prop}The Heisenberg double $H(O,U)$ is isomorphic to the algebra $D_G$ of algebraic differential operators on the group $G$.
\end{prop}
\end{example}

We recall the following well-known and important result:

\begin{thm}\label{faithful} Suppose that $U$ has enough finite dimensional modules.  Then $F$ is a faithful $H(F,U)$-module.
\end{thm}

\begin{proof} Let $\phi,\psi: F\to F$ be linear operators. Recall that the convolution 
$\phi*\psi: F \to F$ is defined as $(\phi*\psi)(f)=\phi(f_1)\psi(f_2)$. Then the action of 
$\phi:=\sum f_i\otimes h_i\in H(F,U)$ on $F$ may be written as $\rho(\phi)=\phi*1$
(where we view $\phi$ as an element of $\End_\CC(F)$ of finite rank, using that 
$U\subset F^*$) It is well known that the operator $\phi\mapsto \phi*1$ on 
$\End_\CC(F)$ is invertible for any Hopf algebra $F$, the inverse being $\phi\mapsto\phi*S$, 
where $S$ is the antipode. Thus $\phi*1=0$ implies $\phi=0$, and we are done. 
\end{proof}

$F$ is sometimes called the basic representation, since when $H(F,U)$ is differential operators on an algebraic group, $F$ is just functions on the group.
\subsection{Quantum $D$-modules on $U$}

In this section, we want to generalize the construction of differential operators to the non-commutative setting.  Our motivating example will be $U=U_t(\gg)$, but we will make the construction for an arbitrary quasi-triangular Hopf algebra $U$.  Key to the present construction is the isomorphism $\Xi$ constructed in \cite{VV}, which relates two potentially different notions of differential operators, one in terms of $F$, the other in terms of $A$.  Indeed, it was Propositions 1.8.1 and 1.8.2 of \cite{VV} which first alerted the author to the relevance of the reflection equation algebra $A$ to the present work.

Recall the reflection equation algebra $A$ constructed previously.  As $A$ is a $U^{[2]}$-equivariant algebra, we can construct the algebra $A\rtimes U^{[2]}$.  

\begin{defn}\cite{VV}\label{DUdefn} The algebra $D_U$ of differential operators on $U$ is the subalgebra $A \ot U \ot 1$ of $A\rtimes U^{[2]}$.  A quantum $D$-module for $U$ is a module over the algebra $D_U$.
\end{defn}

\begin{figure}[htb]
\centerline{\includegraphics[height=5cm]{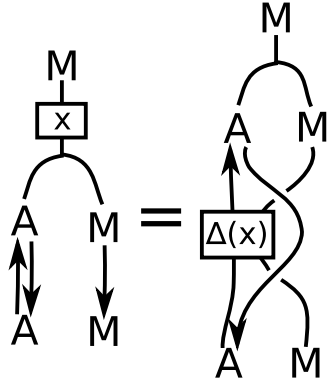}}
\caption{The commutation relations for differential operators.}\label{commreln}
\end{figure}

\begin{rem}\label{dmoduleaction} $M$ is thus both an $A$-module and a $U$-module, such that, for $u\in U, a\in A, m\in M$, we have $x(am)=\sum_{j,k}\mu_M(((x_1\ot S(r_j^+)r_k^+)\rhd a) \ot r_j^-x_2r_k^-m)$.  This commutation relation is depicted graphically in figure \ref{commreln}.
\end{rem}

\begin{rem}\label{co-op}  Typically (e.g. \cite{KrBa}), a quantum $D$-module is defined as a module over the Heisenberg double $H(F,U)$. In \cite{VV}, the authors construct an isomorphism \mbox{$\Xi: D_U\cong H(F,U)$}.  Thus, the present notion of a quantum $D$-module agrees with the usual notion.\footnote{In fact, it will turn out that we could assume slightly less in our constructions:  namely, we could take $M$ to be a module over the subalgebra $A\ot U'\subset D_U$ instead of the full algebra $D_U$, as the constructions of Section 4 only use the action of this subalgebra. However, in this article we will ignore the distinction.}
\end{rem}

\subsection{Three actions of $U$ on $D$-modules}
In this section, we will recreate the left, right, and adjoint actions of vector fields on a classical $D_G$-module, in the non-commutative context.  First, we recall the classical setup.  Let $G$ be an algebraic group, and let $U=U(\gg)$ as before.  Then the left, right and adjoint action of vector fields induce algebra maps $\partial_\rhd:U\to D_U$, $\partial_\lhd: U\to D_U$, and $\ad:U \to D_U$.  We have that $\ad(h)=\partial_\lhd(h_1)\partial_\rhd(h_2)=\mu\circ(\partial_\lhd\ot \partial_\rhd)\Delta(h)$.  The assertion that $\ad$ is a homomorphism relies on the fact that the images of $\partial_\lhd$ and $\partial_\rhd$ centralize one another in $D_U$, since the image of $\partial_\lhd$ is right invariant vector fields, and the image of $\partial_\rhd$ is left invariant vector fields.  A more concise way of saying this is that there is a homomorphism $\phi:U^e\to D_U$, and that the left, right, and adjoint homomorphisms are given by pre-composing with the left, right, and adjoint maps $U\to U^{e}$.  Note that because we act on functions, the left action corresponds to right translation in the group, and vice versa.

In the non-commutative situation, again following \cite{VV}, we define $\partial_\lhd$ and $\partial_\rhd$ as follows.  We have the isomorphism $\Xi:D_U\cong H(F,U)$.  We have the inclusion $\partial_\lhd$ as inclusion into the $U^o$ factor.

\begin{defn}  The adjoint action of $U$ on itself is given by $\ad: x\ot y\mapsto x_1 y S(x_2)$.\end{defn}

\begin{defn} Denote by $U'$ the sub-algebra of $x\in U$ s.t. $\ad(U)x$ is finite dimensional.
\end{defn}
\begin{rem} It is straightforward to check that $U'$ is indeed a subalgebra, and that $\Delta(U')\subset U\ot U'$.
\end{rem}

The adjoint action $U\ot U'\to U'$ of $U$ on itself yields the co-adjoint map $\partial_\rhd: U'\to F\otimes U'\subset H(F,U)$.

\begin{defn}
Let $U^{[2]'}$ denote the subalgebra $U\ot U'$ of $U^{[2]}$.  We define the homomorphism $\partial_2: U^{[2]'}\to H(F,U^o)$ by $x\ot y \mapsto \partial_\lhd(x)\partial_\rhd(y)$.  Abusing notation, we denote also by $\partial_2: U^{[2]'}\to D_U$ the map $x\ot y \mapsto \Xi^{-1}(\partial_\lhd(x)\partial_\rhd(y))$.
\end{defn}

\begin{rem} \label{omegacommutes} Let $z\in U$ be a central element.  $z$ is thus ad-invariant, and so co-ad invariant.   Thus we have that $\partial_\lhd(z)=1\ot z = \partial_\rhd(z)$, in other words the left and right actions of $z$ agree whenever $z$ is central in $U$.  We will apply this observation to the ribbon element later in the construction.
\end{rem}

\section{Statement of Results}
Our main result is that the data of a ribbon Hopf algebra $U$, a f.d. $U$-module $V$ and a quantum $D_U$-module $M$ together yield representations of the elliptic braid group on $n$ strands, for any $n$.  As $D_{U_t}$ is a flat (in fact, trivial) deformation of $D_G$, this provides a rich source of such representations.  Taking $M$ to be the ``basic" $D_U$-module $A$, we recover the construction of Lyubashenko and Majid \cite{Ly}, \cite{LyMa}.  Alternatively, taking a quasi-classical limit as $t\to1$ for the quantum group $U_t(sl_N)$, we recover the geometric constructions from \cite{CEE}.  Thus our results are really an interpolation of those two papers.  In a future paper, we hope to elaborate the case $U=U_t(sl_N)$ in greater detail, and extend some representation-theoretic results of \cite{CEE}.

To state the main theorem, we need to introduce some further notation.  Let $Z$ be a $U$-module.  We denote by $Z^{inv}$ the vector space
$$Z^{inv}=Hom_U(1,Z)=\{w\in Z \textrm{ s.t. } xw=\eps(x)w\forall x\in U\}.$$

Our main result is the following:
\begin{thm}\label{mainresult} Let $U$ be a ribbon Hopf algebra.  Let $n\in \NN$, $V$ a f.d. $U$-module, and $M$ a $D_U$-module. On the vector space 
$$W=(V^{\ot n}\ot M)^{inv}$$
of invariants w.r.t to the adjoint action on $M$, we have an action of the elliptic braid group $B_n^{Ell}$, which defines a functor $F_{n,V}: D_U-mod\to Rep(B_n^{Ell})$.
\end{thm}
We will provide the construction in Section \ref{construction}, from which will follow two easy corollaries:
\begin{corollary} Let $\nu\in U$ be the ribbon element.  Suppose that $V$ is irreducible, and that $\nu|_V=c_V \id_V$.  Then, $\tilde{Y}|_V=\tilde{X}|_V=c_V^n \id_W$.
\end{corollary}
\begin{corollary}\label{DAHA} Suppose the braiding on $V$ satisfies the Hecke relation $$(\sigma_{VV}-q^{-1}t)(\sigma_{VV}+q^{-1}t^{-1}).$$  Then the action of $B_n^{Ell}$ descends to an action of the double affine Hecke algebra (DAHA) $\HH_n(q,t)$, and we have a functor $F_{n,V}:D_U-mod\to Rep(\HH_n(q,t))$.
\end{corollary}
Finally, we consider the example $U=U_t(sl_N)$, and $V$ is the defining representation, of highest weight $(1,0,\ldots,0)$.  It is well-known that $V$ is Hecke with parameters $q,t=q^{nk}, k=N/n$.  In \cite{CEE}, the authors considered $D$-modules on (the classical group) $SL_N$; given a $D$-module $M$, they constructed a representation of the trigonometric Cherednik algebra $\HH^{deg}_n(k)$ of type $A_{n-1}$ on the space $(V^{\ot n}\ot M)^{inv}$, with parameter $k=N/n$.  As $\HH^{deg}_n(k)$ is the quasi-classical limit of the DAHA as the parameter $t\to 1$, we can ask whether our construction agrees with theirs in the quasi-classical limit.  Indeed, we have the following
\begin{thm} In the quasi-classical limit as $t\to 1$ the construction in Proposition \ref{DAHA} recovers the $\HH_n^{deg}(k)$-representations constructed in \cite{CEE}.
\end{thm}

If we forget the $A$-action, and consider $M$ only as a $U^e$-module, i.e. a $U$-bimodule, we can still define the operators $Y_i,T_j$.  In fact, we have the following result which follows from the proof of Theorem \ref{mainresult}

\begin{corollary} Let $M$ be a $U^e$-module.  The operators $Y_i, T_j$ define an action of the affine Hecke algebra on $W$, whose quasi-classical limit is the Arakawa-Suzuki construction from \cite{AS}.
\end{corollary}

\section{The Construction}\label{construction}
This entire section is devoted to a constructive proof of Theorem \ref{mainresult}.  We will build the representation of $B_n^{Ell}$ by constructing first the braid group representation, then the action of the algebra $\CC[Y]$, then the algebra $\CC[X]$, and finally checking the commutation relations between them.

Fix $U$, a ribbon Hopf algebra with ribbon element $\nu$, and let $A$ and $D_U$ be as in Definitions \ref{Adefn} and \ref{DUdefn}.  Let $M$ be a $D_U$-module, let $V$ be a $U$-module, and consider the vector space $$W=(V^{\ot n}\ot M)^{inv}.$$  This means that we take invariants with respect to the usual action on the $V^{\ot n}$ factor, and with respect to the adjoint action on the $M$ factor.  It is important to note here that the adjoint action is only defined for elements of $U'$, the locally finite part of $U$.  So we will only be able to apply ad-invariance for such elements.

We index the tensor factors in $W$ from right to left, starting with $M$ at index $0$.  We will use the symbol $\pi_V$ when we want to explicitly emphasize the action on $V$.  

The action of the braid group is given by $T_i=\sigma_{i+1,i}$ acting on the $V^{\ot n}$ factor.  We define the invertible element $Y_1=\sigma_{M,V}\circ\sigma_{V,M}=(\pi_V\ot\partial_\lhd)(R_{01} R_{10})$ acting on the right two factors.  We define operators $Y_i$, for $i=1,\ldots,n-1$ by $Y_{i+1}=T_iY_iT_i$.  It follows from the QYBE that the $Y_i$'s commute pairwise.  Recall that we denote by $\tilde{Y}$ the product \mbox{$\tilde{Y}=\prod_iY_i$}.

\begin{prop}\label{ribboncor} $$\tilde{Y}=(\pi_V\ot\cdots\ot\pi_V\ot\partial_\lhd)\Big(\Delta^{(n)}(\nu^{-1})\underbrace{(\nu\ot\cdots\ot\nu)}_{n+1}\Big).$$\end{prop}
\begin{proof} This expression for the iterated coproduct  of $\nu$ in terms of the double-braidings is essentially its defining property.
\end{proof}

\begin{figure}[htb]
\centerline{\includegraphics[height=3cm]{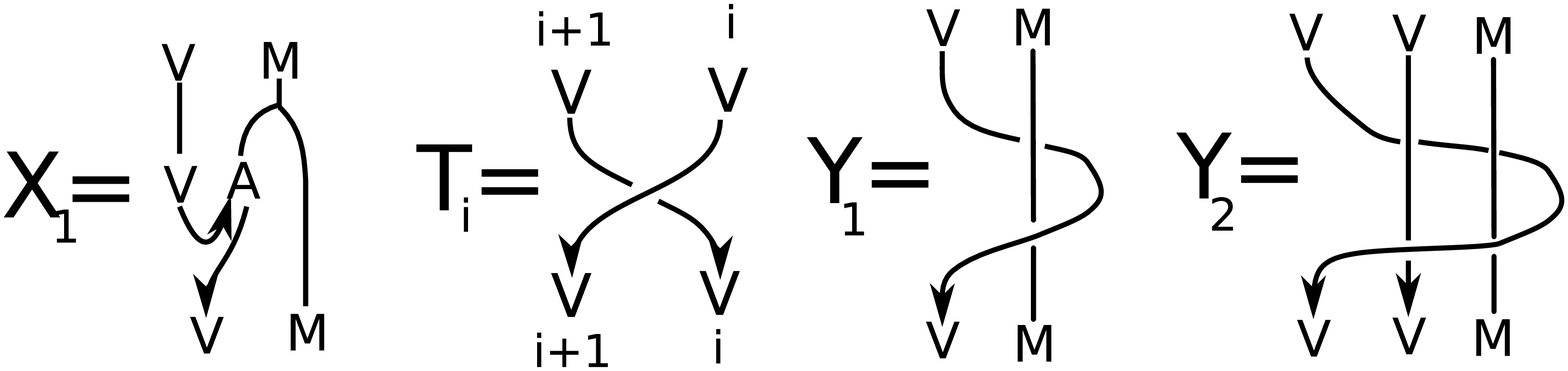}}
\caption{The operators $X_1$, $T_i$, $Y_1$, $Y_2$.}\label{XYTDefn}
\end{figure}

Recall that since $V$ is an $A$-comodule, we have the operators $L\in A\ot End_\CC(V)$ (see Prop \ref{refleqn}).  We can define $X_1=L_{01}$.  Since $\Delta_0(L_{01})=L_{02}L_{12}$, we have that $S_0(L_{01})$ gives an inverse $X_1^{-1}$ (Here $S_0$ means that we apply the antipode in the $A$ component).  We set $X_{i+1}=T_iX_iT_i$, for $i=1,\ldots,n-1$.  The reflection equation implies that the $X_i$ commute pairwise.  The invariant description for $X_1$ is depicted in figure \ref{XYTDefn}, and it is perhaps more enlightening.  We first comultiply $V$ as an $A$-module, then multiply the extra $A$-factor into $M$.

\begin{lemma} $X_i,Y_i, T_i$ preserve the subspace of ad-invariants.\end{lemma}
\begin{proof}
Note that $Y_i$ and $T_i$ are endomorphisms of the $U$-module $V^{\ot n}\ot M$.  Thus, if $v\in (V^{\ot n}\ot M)^{inv}$, we have that $x T_i v=T_i x v= \eps(x) T_iv$, and likewise for $Y_i$.

Only $X_i$ requires some computation.  Since $X_{i+1}=T_iX_iT_i$, we only need to check for $X_1$.  Since $D_U$ is a module-algebra under the adjoint action, it follows that for any $x\in U$, we have
\begin{eqnarray*}
(ad x) X_1(\sum_iv^i_n\ot\cdots\ot v^i_1\ot m^i) &=& (ad x)\mu_{21,0}(\sum_iv^i_n\ot\cdots\ot e^i_k\ot e^{i*}_k\ot v^i_1\ot m^i)\\
&=& \mu_{21,0}(ad x)(\sum_iv_n^i\ot\cdots\ot e^i_k\ot e^{i*}_k\ot v^i_1\ot m^i)\\
&=& \eps(x) \mu_{21,0} (\sum_iv_n^i\ot\cdots\ot e^i_k\ot e^{i*}_k\ot v^i_1\ot m^i)\\
&=&\eps(x) X_1(\sum_iv_n^i\ot\cdots\ot v^i_1\ot m^i)
\end{eqnarray*}
We used that $U$ acts trivially on $e_k\ot e_k^*$, as it is just the image of $1\in\CC$ under the co-evaluation, and that the multiplication in $D_U$ is covariant for the adjoint action of $U$.  
\end{proof}

\begin{prop}$X_1Y_2=Y_2X_1T_1^2$.\label{mainformula}\end{prop}
\begin{proof}
We use remark \ref{dmoduleaction} to express both sides as a braid diagram in figure \ref{braidproof}, at which point the equality is simply an identity of tangles.
\begin{figure}[htb]
\centerline{\includegraphics[height=6cm]{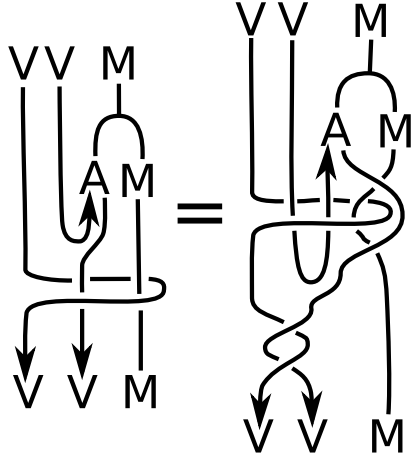}}
\caption{Proof of Proposition \ref{mainformula}}\label{braidproof}
\end{figure}
\end{proof}
This relation appears in a different form in \cite{VV}, where it is computed using the Fourier transform isomorphism: $\mathcal{F}:A\to U'$. Recall that $U'$ is the sub-algebra of $U$ which is locally finite under the adjoint representation.  In our case, $\mathcal{F}$ applied to the zeroeth factor maps the action of $X_1$ onto that of $Y_1$, and thus lines up with the Fourier transform of the torus.  Expressing $X_1$ in terms of the co-product of $A$ allows us to give a diagrammatical proof instead.

\begin{prop}\label{tildeY} $\tilde{Y}X_j= X_j\tilde{Y}$\end{prop}
\begin{proof}
Actually, we'll show that $\tilde{Y}$ acts by $(\nu\ot \ldots \ot \nu)|_{V^{\ot n}}$, so that it clearly commutes with $X_j$.  We have  by corollary \ref{ribboncor} that
\begin{eqnarray*}
\tilde{Y}(\sum_iv_n^i\ot\cdots\ot v_1^i\ot m^i)&=&\Delta^{(n)}(\nu^{-1})(\sum_i\nu v_n^i\ot\cdots\ot \nu v_1^i\ot \partial_\lhd(\nu) m^i)\\&=&\sum_i\nu v_n^i\ot\cdots\ot \nu v_1^i\ot \partial_\lhd(\nu)\partial_\rhd(S^{-1}(\nu^{-1})) m^i, \textrm{ by ad invariance,}\\
&=&\sum_i\nu v_n^i\ot\cdots\ot \nu v_1^i\ot m^i,
\end{eqnarray*}
by remark \ref{omegacommutes}, and the fact that $S^{-1}(\nu)=\nu$.
Notice that this is the only place in the argument where we use ad-invariance.   Because $\nu$ is central, it is in $U'$.  Also it is easy to see that $\Delta^{(k)}(U')\subset U^{\ot n}\ot U'$.   Together, these justify shifting the ribbon element to the zeroeth component.
\end{proof}

We have confirmed the necessary relations on the operators $X_i,Y_i,T_j$, which concludes the proof of Theorem \ref{mainresult}.

\section{A relation for $\tilde{X}$ and $\tilde{Y}$}
As an immediate corollary of the proof of Proposition \ref{tildeY}, we have the following
\begin{corollary} Suppose that $V$ is irreducible, and that $\nu|_V=c_V\id$.  Then $\tilde{Y}=c_V^n \id_W$.
\end{corollary}

It turns out that $\tilde{X}:=\prod_iX_i$ acts by the same scalar.  We prove this in two steps.  Recall our standing assumption that $U$ has enough finite dimensional representations, in the sense of Definition \ref{nondeg}.

\begin{prop} $\tilde{X}=(\tilde{T}\ot1\ot1)\circ\mu_{A,M}\circ \Delta_{V^{\ot n}}$, where $\tilde{T}:V^{\ot n}\to V^{\ot n}$ is the double braiding,
$\tilde{T}=\prod_{i=1}^{n-1}\prod_{j\leq i}T_j$, depicted in figure \ref{tildeT}.

\begin{figure}[htb]
\centerline{\includegraphics[height=3cm]{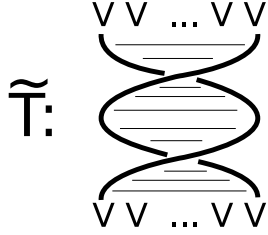}}
\caption{$\tilde{T}$ braids all $n$ copies of $V$ twice around one another.}\label{tildeT}
\end{figure}
\end{prop}
\begin{proof}  This is a direct computation, of which the reflection equation (\ref{refleqn}) is the case \mbox{$n=2$}.
\end{proof}
 
\begin{figure}[htb]
\centerline{\includegraphics[height=6.5cm]{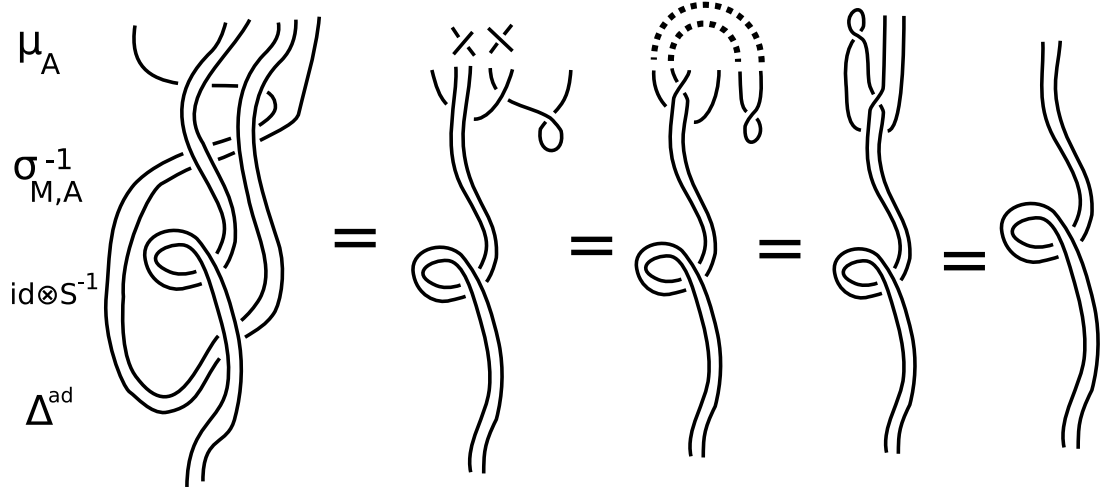}}
\caption{Computation of $\tilde{X}$.  The second equality applies Proposition \ref{auto} to the braiding, as depicted above the diagram.  The third equality applies Proposition \ref{auto} to the braiding composed with the coevaluation, as indicated by the dotted lines.  The first and fourth equalities are already clear at the level of braiding diagrams.\label{Xreln}}
\end{figure}

\begin{prop}\label{tildeX} Suppose that $V$ is irreducible. Then $\tilde{X}=c_V^n$, where $c_V$ is the same as from Proposition \ref{tildeY}.\end{prop}
\begin{proof}  We first consider the case $n=1$.  Since we consider $U$-invariants, this is the same as $A$-coinvariants, and so we have that
\begin{eqnarray}
\Delta_V &=& \sigma^{-1}_{M,A}\circ(\id\ot\id\ot S^{-1})\circ \Delta^{ad}_A, \textrm{ and thus,} \nonumber\\
\tilde{X} &=& \mu_A\circ\sigma^{-1}_{M,A}\circ(\id\ot\id\ot S^{-1})\circ \Delta^{ad}_A,\label{bigterm}
\end{eqnarray}
where $\Delta^{ad}_A$ is the co-adjoint action of $A$ on itself.  In figure \ref{Xreln}, we compute that equation \ref{bigterm} is equal $\ad(\nu)$ acting on $M$.  For concreteness, we work in $M=A$, the basic representation of $D_U$, which we may do because this representation is faithful, by Theorem \ref{faithful}, and our assumption that $U$ has enough finite dimensional representations.

By invariance, $\ad(\nu)$ acting on $M$ is the same as $S(\nu)=\nu$ acting on $V$, as desired.  When $n\geq 1$, we apply this proof to $V^{\ot n}$ instead, to conclude that
\begin{eqnarray*}
\tilde{X}&=&(\tilde{T}\ot1\ot1)\circ\mu_{A,M}\circ \Delta_{V^{\ot n}}\\
&=& \tilde{T}\circ \Delta^{(n-1)}(\nu)\ot1\ot1, \textrm{ by the $n=1$ case above.}\\
&=& (\nu\ot\cdots\ot\nu\ot1\ot1)|_{(V^{\ot n}\ot A)^{inv}}\\
&=& c_V^n \id
\end{eqnarray*}
\end{proof}

\begin{corollary} Let $B_n^{ell,0}$ be the quotient group of $B_n^{ell}$ by the relations $\tilde X=\tilde Y=1$, and $\HH_n^0(q,t)$ be the quotient  of $\HH_n(q,t)$ by the same relations.  Setting $X_i'=X_ic_V^{-1}$, $Y_i'=Y_ic_V^{-1}$, we  obtain an action of the group $B_n^{ell,0}$ on $W$, which descends 
to a representation of $\HH_n^0(q,t)$ when $V$ is Hecke. 
\end{corollary}

\section{DAHA case and its degeneration}
In this section, we consider in more detail the case $U=U_t(sl_N)$.  We recall the degeneration process from the DAHA to trigonometric Cherednik algebra, and we show that our construction degenerates to the construction of \cite{CEE} (with trivial modifications) in the quasi-classical limit.

\subsection{A representation of $\HH_n(q,t)$}
\begin{prop}\label{hecke}  Let $V$ be the defining representation (of highest weight $(1,0,\dots,0)$) for $U_t(sl_N)$, and let $M$ be a $D_U$-module.  Then the operators $X_i$, $Y_i$, and $T_i$ of the preceding section define a representation of $\HH_n(q,t)$ with parameters $q, t=q^{nk}, k=N/n$ on the space $W=(V^{\ot n}\ot M)^{inv}$, of invariants with respect to the adjoint action.\end{prop}
\begin{proof}
It is well known that the operators $\sigma_{V,V}$ satisfy the Hecke relation $$(\sigma_{VV}-t^{1-1/N})(\sigma_{VV}+t^{-1-1/N}),$$ which is the only additional relation on $\HH_n(q,t)$ when $t=q^{nk},k=N/n$.
\end{proof}

\subsection{The trigonometric construction}

In \cite{CEE}, the authors constructed a representation of the trigonometric Cherednik algebra $\HH_n^{deg}(k)$ on the space $(V^{\ot n}\ot M)^{inv}$ of invariants with respect to the adjoint action, where $V$ is the vector representation, and $M$ is a D-module on $G$.  The operators $s_{ij}$ were defined by the usual symmetric group action on $V^{\ot n}$, while the $X_i, y_j$ were defined by:
\begin{eqnarray*}
X_i &=& \sum_{r,s} (E^r_s\ot A^s_r)_{i,0}\\
y_j &=& k (\sum_p (b_p\ot L_{b_p})_{j,0} + \sum_{i<j}s_{ij}),
\end{eqnarray*}
where $A^s_r$ denotes multiplication by the function $A^s_r$ in the coordinate algebra, $\{b_p\}$ are a orthonormal basis with respect to the trace form on $sl_N$, and $L_{b_p}$ denotes the action by left translation along the vector field given by $b_p$.

\subsection{Trigonometric degeneration of the DAHA}
We recall the process of degeneration from the DAHA to the trigonometric Cherednik algebra.  Let $\hbar$ denote a formal variable, let $k\in\CC$, and let $q,t\in \CC[[\hbar]]$ denote the power series
\begin{eqnarray}
q=e^{\hbar}, t=e^{nk\hbar}.\label{qt}
\end{eqnarray}
\begin{defn}$\widetilde{\HH}_n$ is the $\CC$-algebra freely generated by
$X_i^{\pm 1},y_i,s_j$, for \mbox{$i=1,\ldots, n$}, and  \mbox{$j=1,\ldots, n-1$}.  Let $Y_i, T_j\in\widetilde{\HH}_n[[\hbar]]$ denote the power series
$$Y_i=e^{\hbar y_i},\;\;\; T_j=q^{-1}s_je^{\hbar ks_j}.$$
\end{defn}
Let $\tilde{I}$ denote the closed ideal in $\widetilde{\HH}_n[[\hbar]]$ generated by the DAHA relations from Definitions \ref{BnEllRelns} and \ref{DAHAdefn}, and let $I$ denote the saturation of $\tilde{I}$ with respect to $\hbar$, 
$$I = \{x\in \widetilde{\HH}_n[[\hbar]] | \hbar^mx\in \tilde{I} \textrm{ for some $m\geq 0$}\}.$$
\begin{defn} 
$\widehat{\HH}_n(q,t)$ is the quotient:
$$\widehat{\HH}_n(q,t)=\widetilde{\HH}_n[[\hbar]]/I.$$
\end{defn}

\begin{defn} $\HH^f_n(q,t)$ is the complete $\CC[[\hbar]]$-subalgebra of $\widehat{\HH}_n(q,t)$ generated by $X_i, Y_i,$  \mbox{$i=1,\ldots,n$}, and $T_j, j=1\ldots, n-1$.  It is a formal version of $\HH_n(q,t)$ from Definition \ref{DAHAdefn}.
\end{defn}

\begin{prop} \label{extend} Let $V$ be a representation of $\HH^f_n(q,t)$ which is flat as a $\CC[[\hbar]]$-module, and $\rho: \HH^f_n(q,t)\to End(V)$ the corresponding map.  Then $\rho$ extends to a representation of $\widehat{\HH}_n(q,t)$ if and only if $\rho(Y_i) = 1 \mod \hbar, \forall i$.  In this case the extension is unique.
\end{prop}
\begin{proof}  ($\Rightarrow$) is obvious.  To show ($\Leftarrow$), we need to define the action of $y_i$ and $s_i$ on $V$.  First we observe that since $s_i^2=1$, the identity $T_i=q^{-1}s_ie^{\hbar k s_i}$ can simply be solved for $s_i$:
$$ s_i= \frac{qT_i-\sinh(\hbar k)}{\cosh(\hbar k)}.$$
We can also solve the identity $Y_i=e^{\hbar y_i}$ for $y_i$:
$$y_i=\frac{1}{\hbar} \log(1-(1-Y_i)),$$
where the RHS is a well-defined power series in $\hbar$ because we assumed $1-Y_i$ was divisible by $\hbar$.  Uniqueness follows because we have explicitly solved for the $y_i$'s and $s_i$'s.
\end{proof}

\subsection{Trigonometric degeneration of the construction of section \ref{construction}}
In this section, we again consider $t,q$ as elements of $\CC[[\hbar]]$, as in equation (\ref{qt}).
We have the following well-known proposition
\begin{prop} Let $U_1=U_t/\hbar U_t$.  Then $U_1 \cong U(sl_N)$, the classical enveloping algebra.  Furthermore, $U_t$ is a flat deformation of $U_1$.
\end{prop}
\begin{corollary} Let $D_t=D_{U_t}$ and $D_1=D_t/\hbar D_t$.  Then $D_1\cong D_{U(sl_N)}$, and $D_t$ is a flat  (in fact, trivial) deformation of $D_1$.\end{corollary}

\begin{proof}  Indeed, we have that the Hochschild cohomology of $D_G$ is determined by the singular cohomology of $G$: $HH^i(D_G,D_G)=H^i(G,\CC)$. Thus, when $G$ is semi-simple, we have in particular that $HH^2(D_G,D_G)=H^2(G,\CC)=0$, so that $D_G$ admits no non-trivial formal deformations (see e.g, \cite{E}).
\end{proof}

\begin{defn}  Let $C$ denote the Casimir element $C = \sum_p b_pb_p \in U(sl_N)$, where $\{b_p\}$ form an orthonormal basis with respect to the trace form. Let  $\Omega$ denote the canonical 2-tensor in $U(sl_N)\ot U(sl_N)$, $$\Omega = \frac{\Delta(C)-C\ot 1 - 1\ot C}{2} = \sum_p b_p \ot b_p.$$ \end{defn}

\begin{defn} Let $r$ denote the classical $r$-matrix for $U(sl_N)$, so $\Omega = r_{10} + r_{01}$
\end{defn}

\begin{prop} We have the following relations between $r$,$R$,$R_{10}R_{01}$, and $\Omega$.
\begin{eqnarray*}
R &=& t^r = 1 + k\hbar r \mod \hbar^2,\\
R_{10}R_{01}&=& t^{r_{10}}t^{r_{01}} =1\ot 1 + k\hbar\Omega \mod \hbar^2.\\
\end{eqnarray*}
\end{prop}

\begin{prop} $\Omega$ acts as $\sigma_{VV}-1/N$ on $V\ot V$, the tensor square of the defining representation for $U(sl_N)$.
\end{prop}
\begin{proof}  First we compute the canonical 2-tensor $\widetilde{\Omega}$ for $U(gl_N)$.  Instead of an orthonormal basis $\{b_p\}$, we can choose the basis $E^i_j$ and dual basis $E^j_i$, and so
$$\widetilde{\Omega}= \sum_{i,j} E^i_j\ot E^j_i,$$
which is exactly the flip.  Since $\Omega=\widetilde{\Omega}-1/N(id_V\ot id_V)$, the claim follows.
\end{proof}

\begin{prop} The $\HH_n(q,t)$-representation $W$ of Proposition \ref{hecke} extends to a representation of $\widehat{\HH}_n(q,t)$.\end{prop}
\begin{proof} $Y_i$ is expressed as a product of $R$-matrices, each of which is congruent to $1 \mod \hbar$.  Thus the condition of Proposition \ref{extend} is satisfied.
\end{proof}

\begin{prop} Let $W_1=W/\hbar W=(V_1^{\ot n}\ot M_1)^{inv}$.  The operators $s_i$ act on $W_1$ as the flip $s_{i,i+1}$ of tensor factors.  The operators $X_i$ act as
$$X_i=\sum_{k,l} (E^k_l)_i\ot A^l_k.$$
\end{prop}
\begin{proof} Straightforward computation.\end{proof}

\begin{prop}\label{yis} The operators $y_i$ act as $$y_i= k(\Omega_{i,0} + \sum_{j<i}s_{ij} - \frac{i-1}{N})$$\end{prop}
\begin{proof}
To see the claim for $y_1$, we note that $$Y_1 = (R_{10}R_{01})_{1,0}=1\ot 1 + \hbar\Omega_{1,0} \mod \hbar^2.$$ We proceed by induction:
\begin{eqnarray*}
Y_i &=& T_iY_{i-1}T_i\\
&=& s_i(1+k\hbar r_{i,i-1})(1+k\hbar(\Omega_{i-1,0} +\sum_{j<i-1}s_{i-1,j} - \frac{i-2}{N}))s_i(1 + k\hbar r_{i,i-1}) \mod \hbar^2\\
&=& 1+ k\hbar\left(\sum_{j<i-1} s_{i,j} +  r_{i,i-1} + r_{i-1,i}  + \Omega_{i,0} - (i-2)/N\right) \mod \hbar^2\\
&=& 1 + k\hbar\left(\sum_{j<i-1}s_{i,j} + s_{i,i-1} -\frac{1}{N} + \Omega_{i,0} - (i-2)/N\right) \mod \hbar^2\\
&=& 1 + k\hbar\left(\sum_{j<i}s_{i,j} + \Omega_{i,0} -\frac{i-1}{N}\right) \mod \hbar^2\\
\end{eqnarray*}
\end{proof}
Comparing these with the operators of \cite{CEE}, we see that the quasi-classical limit of the present construction recovers the construction there, up to adding constants.\footnote{Also, our operators $L_b$ act by right instead of left translations, though this is just convention.}

\end{document}